\documentclass[11pt]{article}

\usepackage{algorithm,algorithmic}
\usepackage{comment}
\usepackage{amsfonts,amsmath,amssymb,amsthm}
\usepackage{verbatim}
\usepackage{color}
\usepackage{caption,subcaption}
\usepackage{eepic,epic,epsfig,epstopdf}
\usepackage[]{graphics}

\title{Multilevel Markov Chain Monte Carlo Method for High-Contrast Single-Phase Flow Problems}
\author{Yalchin Efendiev\thanks{Department of Mathematics, Texas A\&M University, College Station, TX 77843-3368, USA.}\and Bangti Jin\thanks{Department of Mathematics, University of California, Riverside, 900 University Ave.
Riverside, CA 92521, USA}\and Michael Presho\footnotemark[1]\and Xiaosi Tan\footnotemark[1]}

\date{\today}

\theoremstyle{break}

\newtheorem{remark}{Remark}[section]

\def\EE{{\mathbb E} }

\newcommand{\bigE}{\mathcal{E}}
\newcommand{\bigD}{\mathcal{D}}

\begin{document}
\maketitle

\begin{abstract}
In this paper we propose a general framework for the uncertainty quantification of quantities of
interest for high-contrast single-phase flow problems. It is based on the generalized multiscale
finite element method (GMsFEM) and multilevel Monte Carlo (MLMC) methods. The former provides a
hierarchy of approximations of different resolution, whereas the latter gives an efficient way to
estimate quantities of interest using samples on different levels. The number of basis functions
in the online GMsFEM stage can be varied to determine the solution resolution and the computational
cost, and to efficiently generate samples at different levels. In particular, it is cheap to generate
samples on coarse grids but with low resolution, and it is expensive to generate samples on fine
grids with high accuracy. By suitably choosing the number of samples at different levels, one can
leverage the expensive computation in larger fine-grid spaces toward smaller coarse-grid spaces, while retaining the
accuracy of the final Monte Carlo estimate. Further, we describe a multilevel Markov chain Monte Carlo method, which sequentially
screens the proposal with different levels of approximations and reduces the number of evaluations
required on fine grids, while combining the samples at different levels to arrive at an accurate estimate.
The framework seamlessly integrates the multiscale features of the GMsFEM with the multilevel feature
of the MLMC methods following the work in \cite{ketelson2013}, and our numerical experiments illustrate its efficiency and accuracy in comparison
with standard Monte Carlo estimates.
\end{abstract}

\section{Introduction}

Uncertainties in the description of reservoir lithofacies, porosity and permeability
are major contributors to the uncertainties in reservoir performance forecasting. The uncertainties can be
reduced by integrating additional data, especially dynamic ones such as pressure or production data,
in subsurface modeling. The incorporation of all available data is essential for the
reliable prediction of subsurface properties. The Bayesian approach provides a principled
framework for combining the prior knowledge with dynamic data in order to make predictions on quantities of
interest. However, it poses significant computational challenges largely due to the fact that exploration
of the posterior distribution requires a large number of forward simulations. High-contrast
flow is a particular example, where the forward model is multiscale in nature and only a limited number of forward simulations
can be carried out before becoming prohibitively expensive. In this paper, we present a framework for
uncertainty quantification of quantities of interest based on the
generalized multiscale finite element method (GMsFEM) and multilevel Monte Carlo (MLMC) methods.
The GMsFEM provides a hierarchy of approximations to the solution, and the MLMC provides an efficient
way to estimate quantities of interest using samples on multiple levels. Therefore, the
framework naturally integrates the multilevel feature of the MLMC with the multiscale nature
of the high-contrast flow problem.

Multiscale methods represent a class of coarse-grid solution techniques that
have garnered much attention over the past two decades (see, e.g., \cite{aarnes04,
apwy07,eh09,hw97,hughes98,jennylt03}). They all hinge on the construction of a coarse
solution space that is spanned by a set of multiscale basis functions. In this
paper we follow the framework of the Multiscale Finite Element Method (MsFEM) \cite{hw97} in
which the basis functions are independently pre-computed, and are obtained through
solving a set of local problems that mimic the global operator, in the hope
of capturing the fine scale behavior of the global system. Then
a global formulation is used to construct a reduced-order solution.
While standard multiscale methods have proven very effective for a variety of applications
\cite{eghe05,eh09,ehg04,gppw11}, recent work has offered a generalized framework
for enriching coarse solution spaces in case of parameter-dependent problems, where the parameter
reflects the uncertainties of the system. Specifically, the GMsFEM is a robust solution
technique in which the standard solution spaces from the MsFEM may be systematically enriched
to further capture the fine behavior of the fine grid solution \cite{bl11,egt11,egw11}. The additional basis functions
are chosen based on localized eigenvalue problems. 

The GMsFEM achieves efficiency via coarse space enrichment, which is split into two stages, 
following an offline-online procedure (see also \cite{barrault,boyaval08,nguyen08,Patera}).
At the first stage of the computation, a larger-dimensional (relative to
the online space) parameter-independent offline space is formed. The offline space accounts
for a suitable range of parameter values that may be used in the online stage, and
constitutes a one-time preprocessing step. The offline space is created by first
generating a set of ``snapshots'' in which a number of localized problems are solved
on each coarse subdomain for a number of parameter values. The offline space is then
obtained through solving localized eigenvalue problems that use averaged parameter
quantities within the space of snapshots. A number of eigenfunctions are
kept in order to form the offline space. At the online stage, we solve analogous
eigenvalue problems using a fixed parameter value within the offline space to form
a reduced-order online space. A notable advantage of the GMsFEM construction is that
flexible coarse space dimension naturally provides a hierarchy of approximations to
be used within the MLMC framework. Further,
we avoid unnecessary large-dimensional eigenvalue computations for each parameter
realization. In particular, the offline stage constitutes a one-time preprocessing
step in which the effect of a suitable range of parameter values is embedded into
the offline space. In turn, the online stage only requires solving much
smaller eigenvalue problems within the offline space, along with the construction of a
conforming or non-conforming basis set.

Multilevel Monte Carlo (MLMC) was first introduced by Heinrich \cite{Heinrich:2001} for high-dimensional
parameter-dependent integrals and was later applied to stochastic ODEs by Giles \cite{Giles:2008b,Giles:2008a},
and PDEs with stochastic coefficients by Schwab et al. \cite{BarthSchwabZollinger:2011} and
Cliffe et al. \cite{CliffeGilesScheichlTeckentrup:2011}. However, it has not been considered in
the ensemble level method context before. The main idea of MLMC methods is
to use a respective number of samples at different levels to compute the expected values of quantities
of interest. In these techniques, more realizations are used at the coarser levels with inexpensive
forward computations, and fewer samples are needed at the finer and more expensive levels due to
the smaller variances. By suitably choosing the number of realizations at each level, one can obtain
a multilevel estimate of the expected values at much reduced computational efforts. It also admits the interpretation as
hierarchical control variates \cite{Speight:2009}. Such hierarchical control variates
are employed in evaluating quantities of interest using the samples from multilevel
distributions.

In this work, we couple the GMsFEM with the MLMC methods to arrive at a general framework
for the uncertainty quantification of the quantities of interest in high-contrast flow problems. Specifically,
we take the dimension of the multiscale space to be the MLMC level, where the accuracy of the
global coarse-grid simulations depends on the dimension of the multiscale coarse space.
The convergence with respect to the coarse space dimension plays a key role in selecting the
number of samples at each level of MLMC. To this end, we take different numbers of online
basis functions to generate the multiscale coarse spaces, running more forward coarse-grid
simulations with the smaller dimensional multiscale
spaces and fewer simulations with larger dimensional multiscale spaces.
By combining these simulation results in a MLMC framework one can achieve better accuracy
at the same cost as the classical Monte Carlo (MC) method. To this end, one needs to assess the convergence
of ensemble level methods  with respect to the coarse space dimension, which can be estimated
based on a small number of a priori computations.

Further, we will consider the use of MLMC jointly with multilevel Markov chain Monte Carlo (MLMCMC) methods following \cite{ketelson2013}. The main
idea of MLMCMC approach is to condition the quantities of interest at one level (e.g., at a finer level)
to that at another level (e.g., at a coarser level). The multiscale model reduction framework
provides the mapping between the levels, and it can be used to estimate the expected value. Specifically,
for each proposal, we run the simulations at different levels
to screen the proposal and accept it conditionally at these levels. In this manner, we obtain
samples from hierarchical posteriors corresponding to our multilevel approximations which can be
used for rapid computations within a MLMC framework.

The rest of the paper is organized as follows. In Section \ref{sec:prelim}, we describe the Bayesian
formulation for the uncertainty quantification of quantities of interest for flow problems, and
in Section \ref{sec:gmsfem}, we describe the two-stage
procedure of the GMsFEM for high-contrast single-phase flow problems. We shall discuss the
offline and online computations in detail. In Section \ref{sec:mlmc} we discuss the idea of
multilevel Monte Carlo methods, and also the crucial issue of complexity analysis. The algorithm
for coupling the GMsFEM with the MLMC is described. Then in Section \ref{sec:mlmcmc}, we
describe a multilevel Markov chain Monte Carlo method for generating samples from
hierarchical posteriori distributions, which can be used in the MLMC framework. A preliminary analysis
of the convergence of the MLMCMC algorithm is also provided. In Section \ref{sec:numer}, we present
numerical examples to illustrate the efficiency of the framework, in comparison with the standard MCMC
estimates. We offer some concluding remarks in Section \ref{sec:concl}.

\section{Preliminaries}
\label{sec:prelim}
Let $D\subset \mathbb{R}^d$ ($d=2,3$) be an open bounded domain, with a boundary $\partial D$.
The model equation for a single-phase, high-contrast flow reads:
\begin{equation}\label{eqn:model}
  -\nabla\cdot (k(x;\mu)\nabla u) = f\quad \mbox{ in } \, \, D
\end{equation}
subject to suitable boundary conditions, where $f$ is the source term and $u$ denotes the pressure
within the medium. Here $k(x;\mu)$ is the heterogeneous spatial permeability field with multiple
scales and high contrast, where $\mu$ represents the dependence on a multidimensional
random parameter, typically resulting from a finite-dimensional noise assumption on the underlying
stochastic process \cite{BabuskaTemponeZouraris:2004}. In practice, one can measure the
observed data $F_{obs}$ (e.g., pressure or production), and then conditions the permeability
field $k$ with respect to the measurements $F_{obs}$ for predicting quantities of interest. Below we recall preliminaries of the
Bayesian formulation (likelihood, prior and posterior) for systematically performing the task.


In this paper, the main objective is to sample the permeability field conditioned on the observed pressure
data $F_{obs}$. The pressure is an integrated response, and the
map from the pressure to the permeability field is not one-to-one. So there may exist many different
permeability realizations that equally reproduce the given pressure data $F_{obs}$. In practice, the
measured pressure data $F_{obs}$ inevitably contains measurement errors. For a given permeability field $k$, we denote
the pressure as $F(k)$, which can be computed by solving the model equation \eqref{eqn:model} on the fine
grids. The computed pressure $F(k)$ will contain also modeling error, which induces an additional source
of errors, apart from the inevitable measurement error. By assuming the combined error as a random variable $\epsilon$
we can write the model as
\begin{equation}
   F_{obs}=F(k)+\epsilon.
\end{equation}
For simplicity, the noise $\epsilon$ will be assumed to follow a normal distribution $\mathcal{N}(0,\sigma_{f}^{2}I)$, i.e., the
likelihood $p(F_{obs}|k)$ is assumed be of the form
\begin{equation*}
  p(F_{obs}|k) \propto e^{-\frac{\|F(k)-F_{obs}\|^2}{2\sigma_{f}^2}}.
\end{equation*}
We will represent the permeability field $k$, which includes facies and interfaces, through the
now classical Karhunen-Lo\'{e}ve expansion (KLE) \cite{Loeve:1977}, which we describe in more detail
in Section \ref{sec:numer}. We let the vector $\theta$ parameterize
the permeability field within facies and $\tau$ parameterize the
velocity in the level set method, respectively. By parameterizing
the interfaces with level sets, the permeability field $k$ is completely determined by $\theta$ and $\tau$.
Our goal is to generate permeability fields $k$ consistent with the observed pressure data $F_{obs}$.
This can be achieved using Bayes' formula which expresses the posterior distribution $\pi(k)$ as
\begin{equation*}
  \begin{aligned}
     \pi(k)= p(k|F_{obs})&\propto p(F_{obs}|k)p(k)\\
        &= p(F_{obs}|k)p((\theta,\tau)')\\
        &= p(F_{obs}|k)p(\theta)p(\tau), 
 \end{aligned}
\end{equation*}
where the last line follows from the standing assumption that the random variables $\theta$ and
$\tau$ are independent. In the expression for the posterior distribution $\pi(k)$, $p(F_{obs}|k)$
is the likelihood function, incorporating the information in the data $F_{obs}$, and $p(\theta)$
and $p(\tau)$ are the priors for the parameters $\theta$ and $\tau$, respectively, encoding prior
knowledge on the permeability fields. In the absence of interfaces, we shall suppress the notation
$\tau$ in the above formula. Further, we may also incorporate other prior information, e.g., that
the permeability field $k$ is known at some spatial locations corresponding to wells.


From a probabilistic point of view, the problem of sampling from the posterior distribution
$\pi(k)$ amounts to conditioning the permeability fields to the pressure data $F_{obs}$
with measurement errors, i.e., the conditional distribution
$p(k|F_{obs})$. 
Generally, this is achieved by Markov chain Monte Carlo (MCMC) methods,
especially the Metropolis-Hastings algorithm.
The main computational effort of the algorithm lies in
evaluating the target distribution $\pi(k)$, which enters the computation through
the acceptance probability. The map between the permeability $k$
and the pressure data $F(k)$ is only defined implicitly by the governing PDE system.
Hence, to evaluate the acceptance probability, one needs to solve
a PDE system on the fine-scale for any given permeability $k$.
Consequently, its straightforward application is very
expensive, which necessitates the development of faster algorithms. In the next section,
we describe the GMsFEM for the efficient forward simulation to provide a hierarchy
of approximations which can be efficiently used for constructing Monte Carlo estimates.

\section{GMsFEM} \label{sec:gmsfem}

\begin{figure}[tb]
  \centering
  \includegraphics[width=0.7 \textwidth]{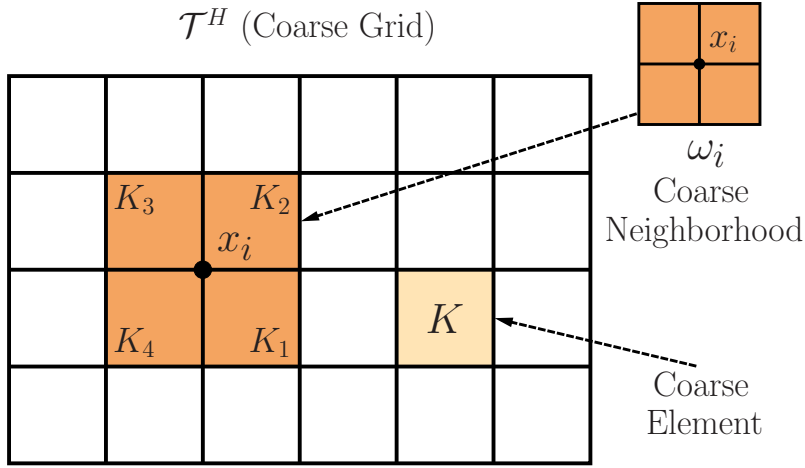}
  \caption{Illustration of a coarse neighborhood and coarse elements}
  \label{schematic}
\end{figure}

To discretize the model equation \eqref{eqn:model}, we first introduce the notion of fine and
coarse grids. Let $\mathcal{T}^H$ be a conforming triangulation of the computational domain $D$ into
finite elements (triangles, quadrilaterals, tetrahedra, etc.). We refer to this partition
as the coarse grid and assume that each coarse subregion is further partitioned into a
connected union of fine grid blocks. The fine grid partition will be denoted by $\mathcal{T}^h$.
We use $\{x_i\}_{i=1}^{N_v}$ (where $N_v$ is the number of coarse nodes) to denote the vertices of
the coarse mesh $\mathcal{T}^H$, and define the neighborhood $\omega_i$ of the node $x_i$ by
\begin{equation} \label{neighborhood}
   \omega_i=\bigcup\{ K_j\in\mathcal{T}^H; ~~~ x_i\in \overline{K}_j\}.
\end{equation}
See Fig.~\ref{schematic} for an illustration of neighborhoods and elements subordinated to the coarse discretization.

Next we describe the offline-online computational procedure for the efficient construction
of GMsFEM coarse spaces. At the offline stage, one generates the snapshot set, and constructs
a low-dimensional offline space by model reduction. At the online
stage, for each input parameter $\mu$, one first computes multiscale basis functions and then
solves for a coarse-grid problem for any force term and boundary condition.
Below we describe the offline and online procedures in more detail.

\subsection{Offline computation}
At the offline stage, we first construct a snapshot space $V_{\text{snap}}^{\omega_i}$ on
each coarse neighborhood $\omega_i$ in the domain (cf. Fig.~\ref{schematic}). The construction
involves solving a set of localized problems for various choices of input parameters. Specifically, we
solve the following eigenvalue problems on each $\omega_i$:
\begin{equation} \label{snaplinalg}
A(\mu_j) \psi_{l,j}^{\omega_i, \text{snap}} = \lambda_{l,j}^{\omega_i, \text{snap}} S(\mu_j) \psi_{l,j}^{\omega_i, \text{snap}}
\quad \text{in} \, \, \, \omega_i,
\end{equation}
where $\{\mu_j\}_{j=1}^J$ is a set of parameter values to be specified. Here we consider only
Neumann boundary conditions, but other boundary conditions are also possible.
The matrices $A(\mu_j)$ and $S(\mu_j)$ in \eqref{snaplinalg} are respectively defined by
\begin{equation}\label{eqn:eigenmatrix}
   \begin{aligned}
       A(\mu_j) & = [a(\mu_j)_{mn}] = \int_{\omega_i} \kappa(x; \mu_j) \nabla \phi_n \cdot \nabla \phi_m dx,\\
       S(\mu_j) &= [s(\mu_j)_{mn}] = \int_{\omega_i} \widetilde{\kappa}(x; \mu_j)  \phi_n \phi_mdx,
   \end{aligned}
\end{equation}
where $\phi_n$ denotes the standard bilinear, fine-scale basis functions and $\widetilde{k}$
will be described below, cf. \eqref{ktilda}. We note that \eqref{snaplinalg} is the discrete
counterpart of the continuous Neumann eigenvalue problem
\begin{equation*}
-\text{div}(\kappa(x, \mu_j) \nabla \psi_{l,j}^{\omega_i, \text{snap}} ) = \lambda_{l,j}^{\omega_i, \text{snap}}
   \widetilde{\kappa}(x; \mu_j) \psi_{l,j}^{\omega_i, \text{snap}} \quad \text{in} \, \, \, \omega_i.
\end{equation*}
For notational simplicity, we omit the superscript $\omega_i$. For each $\omega_i$, we keep the
first $L_i$ eigenfunctions of \eqref{snaplinalg}
corresponding to the lowest eigenvalues to form the snapshot space
\begin{equation*}
  V_{\text{snap}}= \text{span}\{ \psi_{l,j}^{ \text{snap}}:~~1\leq j \leq J, \ 1\leq l \leq L_i \}.
\end{equation*}
We then stack the snapshot functions into a matrix
\begin{equation*}
   R_{\text{snap}} = \left[ \psi_{1}^{\text{snap}}, \ldots, \psi_{M_{\text{snap}}}^{\text{snap}} \right],
\end{equation*}
where $M_{\text{snap}}=J\times L_i$ denotes the total number of snapshots used in the construction.

Next we construct the offline space $V_{\text{off}}^{\omega_i}$, which will be
used to efficiently (and accurately) construct a set of multiscale basis functions
for each $\mu$ value at the online stage. To this end, we perform a dimensionality reduction
of the snapshot space using an auxiliary spectral decomposition. Specifically, we seek a subspace of the snapshot space such that
it can approximate any element of the snapshot space in a suitable sense. 
The analysis in
\cite{egw11} motivates the following eigenvalue problem in the space of snapshots:
\begin{equation} \label{offeig}
A^{\text{off}} \Psi_k^{\text{off}} = \lambda_k^{\text{off}} S^{\text{off}} \Psi_k^{\text{off}},
\end{equation}
where the matrices $A^\mathrm{off}$ and $S^\mathrm{off}$ are defined by
\begin{equation*}
  \begin{aligned}
    A^{\text{off}} &= [a_{mn}^{\text{off}}] = \int_{\omega_i} \overline{\kappa}(x; \mu)
   \nabla \psi_m^{\text{snap}} \cdot \nabla \psi_n^{\text{snap}}dx = R_{\text{snap}}^T \overline{A} R_{\text{snap}},\\
    S^{\text{off}} &= [s_{mn}^{\text{off}}] = \int_{\omega_i } \overline{\widetilde{\kappa}}(x; \mu)\psi_m^{\text{snap}} \psi_n^{\text{snap}}dx = R_{\text{snap}}^T \overline{S} R_{\text{snap}},
 \end{aligned}
\end{equation*}
respectively. Here $\overline{\kappa}(x, \mu) $ and $\overline{\widetilde{\kappa}}(x, \mu)$ are
domain-based, parameter-averaged coefficients, and $\overline{A}$ and $\overline{S}$ denote fine
scale matrices for the averaged coefficients. To generate the offline space, we choose the smallest
$M_{\text{off}}$ eigenvalues to \eqref{offeig}, and take the corresponding eigenvectors in the space of snapshots
by setting $\psi_k^{\text{off}} = \sum_j \Psi_{kj}^{\text{off}} \psi_j^{\text{snap}}$ (for $k=1,\ldots,
M_{\text{off}}$) to form the reduced snapshot space, where $\Psi_{kj}^{\text{off}}$ are the coordinates of
the vector $\Psi_{k}^{\text{off}}$. We then create the offline matrix
\begin{equation*}
R_{\text{off}} = \left[ \psi_{1}^{\text{off}}, \ldots, \psi_{M_{\text{off}}}^{\text{off}} \right]
\end{equation*}
to be used in the online computation.

\begin{remark}
At the offline stage the bilinear forms in \eqref{offeig} are chosen to be parameter-independent,
such that there is no need to construct the offline space for each $\mu$ value.
\end{remark}

\subsection{Online computation}
Next for a given input parameter $\mu$ value, we construct the associated online coarse space
$V^{\omega_i}_{\text{on}}(\mu)$ on each coarse subdomain $\omega_i$.
In principle, we want this to be a low-dimensional subspace of the offline space for computational efficiency.
The online coarse space will be used by the continuous Galerkin finite element method
for solving the original global problem.
In particular, we seek a subspace of the offline space such that it can approximate any element of the offline space in
an appropriate sense. We note that at the online stage, the bilinear forms are \emph{parameter-dependent}.
The analysis in \cite{egw11} motivates the following eigenvalue problem in the offline space:
\begin{equation} \label{oneig}
A^{\text{on}}(\mu) \Psi_k^{\text{on}} = \lambda_k^{\text{on}} S^{\text{on}}(\mu) \Psi_k^{\text{on}},
\end{equation}
where the matrices $A^\mathrm{on}(\mu)$ and $S^\mathrm{on}(\mu)$ are defined by
\begin{equation*}
   \begin{aligned}
    A^{\text{on}}(\mu) &= [a^{\text{on}}(\mu)_{mn}] = \int_{\omega_i} \kappa(x; \mu) \nabla \psi_m^{\text{off}} \cdot \nabla \psi_n^{\text{off}} dx= R_{\text{off}}^T A(\mu) R_{\text{off}},\\
    S^{\text{on}}(\mu) &= [s^{\text{on}}(\mu)_{mn}] = \int_{\omega_i} \widetilde{\kappa}(x; \mu) \psi_m^{\text{off}} \psi_n^{\text{off}} dx = R_{\text{off}}^T S(\mu) R_{\text{off}},
   \end{aligned}
 \end{equation*}
respectively. Note that $\kappa(x; \mu)$ and $\widetilde{\kappa}(x; \mu)$ are now parameter dependent. To
generate the online space, we choose the eigenvectors corresponding to the smallest $M_{\text{on}}$
eigenvalues of \eqref{oneig}, and set
$\psi_k^{\text{on}} = \sum_j \Psi_{kj}^{\text{on}} \psi_j^{\text{off}}$ (for $k=1,\ldots, M_{\text{on}}$),
where $\Psi_{kj}^{\text{on}}$ are the coordinates of the vector $\Psi_{k}^{\text{on}}$.

\begin{remark}[{Adaptivity in the parameter space}]

We note that one can
use adaptivity in the parameter space to avoid computing
the offline space for a large range of parameters and compute the
offline space only for a short range of parameters and update the space.
To demonstrate this concept, we assume that the parameter space
$\Lambda$ can partitioned into a number of smaller parameter spaces
$\Lambda_i$, $\Lambda=\bigcup_i \Lambda_i$, where
$\Lambda_i$ may overlap with each other. Furthermore, the offline
spaces are constructed for each $\Lambda_i$. In the online stage,
depending on the online value of the parameter, we can decide
which offline space to use. This reduces the computational cost
in the online stage. In many applications, e.g., in nonlinear problems,
one may remain in one of $\Lambda_i$'s for many iterations and
thus use the same offline space to construct the online space.
Moreover, one can also adaptively add multiscale basis functions
in the online stage using error estimators. This is a subject of our
future research.
\end{remark}

\subsection{Global coupling mechanism}
To incorporate the online basis functions into a reduced-order global
formulation of the original problem \eqref{eqn:model}, we begin with an initial coarse space
$V^{\text{init}}(\mu) = \text{span}\{ \chi_i \}_{i=1}^{N_v}$ ($N_v$ denotes
the number of coarse nodes). The functions $\chi_i$ are standard multiscale partition
of unity functions defined by
\begin{equation*}
  \begin{aligned}
   -\text{div} \big( \kappa(x; \mu) \, \nabla \chi_i   \big) = 0 \quad K \in \omega_i, \\
    \chi_i = g_i \quad \text{on} \, \, \, \partial K,
  \end{aligned}
\end{equation*}
for each coarse element $K \in \omega_i$, where $g_i$ is a bilinear boundary condition. Next we 
define the summed, pointwise energy $\widetilde{\kappa}$ as, cf. \eqref{eqn:eigenmatrix},
\begin{equation} \label{ktilda}
\widetilde{\kappa} = \kappa \sum_{i=1}^{N_v} H^2 | \nabla \chi_i |^2,
\end{equation}
where $H$ is the coarse mesh size. In order to construct the global coarse grid solution space
we multiply the partition of unity functions $\chi_i$ by the online eigenfunctions $\psi_k^{\omega_i,
\text{on}}$ from the space $V^{\omega_i}_{\text{on}}(\mu)$ to form the basis functions
\begin{equation} \label{finalbasis}
\psi_{i,k} = \chi_i \psi^{\omega_i, \text{on}}_k \quad \text{for} \, \, 1 \leq i \leq N_v \, \, \text{and} \, \,
1 \leq k \leq M_{\text{on}}^{\omega_i},
\end{equation}
where we recall that $M_{\text{on}}^{\omega_i}$ denotes the number of online basis functions kept
for each $\omega_i$. The basis constructed in \eqref{finalbasis} is then used
within a global continuous Galerkin formulation. Now we define the online spectral multiscale space as
\begin{equation}
V_{\text{on}}(\mu) = \text{span}\{  \psi_{i,k}: \, 1 \leq i \leq N_v, \, \,  1 \leq k \leq M^{\omega_i}_{\text{on}}    \},
\end{equation}
and using a single index notation, we write $V_{\text{on}}(\mu) = \text{span}\{  \psi_{i}   \}_{i=1}^{N_c}$
where $N_c$ denotes the total number of basis functions in the coarse scale formulation. Using the online
basis functions, we define the operator matrix $R = [\psi_1, \ldots, \psi_{N_c}]$, where $\psi_i$ represents
the vector of nodal values of each basis function defined on the fine grid. To solve \eqref{eqn:model}
we seek $u(x; \mu) = \sum_i u_i \psi_i(x; \mu) \in V_{\text{on}}$ such that
\begin{equation} \label{coarsevarform}
  \int_D \kappa(x;\mu)\nabla u \cdot \nabla v dx = \int_D f \, v dx\, \, \, \text{for all} \, \, v \in V_{\text{on}}.
\end{equation}
The above equation yields the discrete form
\begin{equation} \label{coarselinalgform}
  A(\mu) u = F,
\end{equation}
where $\displaystyle A(\mu) := [a_{IJ}] = \int_D \kappa(x; \mu) \nabla \psi_I \cdot \nabla \psi_Jdx$ is a
coarse stiffness matrix, $\displaystyle F := [f_I] = \int_D f  \, \psi_Idx$ is the coarse forcing vector,
$P_c$ denotes the vector of unknown pressure values, and $\psi_I$ denotes the coarse basis functions that
span $V_{\text{on}}$. We note that the coarse system may be rewritten using the fine-scale system and the
operator matrix $R$. In particular, we may write $A(\mu) = R^T A^f(\mu) R$ and $F =R^T F^f$, where
\begin{equation*}
   A^f(\mu) := [a_{ij}] = \int_D \kappa(x; \mu) \nabla \phi_i \cdot \nabla \phi_jdx\quad \mbox{and}\quad
   F^f := [f_i] = \int_D f  \, \phi_idx,
\end{equation*}
and $\phi_i$ are the fine-scale bilinear basis functions. Analogously, the operator matrix $R$ may be used
to map coarse scale solutions back to the fine grid.

\section{Multilevel Monte Carlo methods}
\label{sec:mlmc}

As was mentioned earlier, one standard approach for exploring posterior distributions is the
Monte Carlo method, especially Markov chain Monte Carlo (MCMC) methods. Here generating each
sample requires the solution of the forward model, which is unfortunately very expensive for many
practical problems defined by partial differential equations, including high contrast flows. Therefore,
it is imperative to reduce the computational cost of the sampling step. We shall couple the multilevel
Monte Carlo with the multiscale forward solvers to arrive at a general framework for uncertainty
quantification of high-contrast flows.

\subsection{MLMC-GMsFEM framework}\label{subsec:mlmcbg}

The MLMC approach was first introduced by Heinrich in \cite{Heinrich:2001} for finite- and infinite-dimensional
integration. Later on, it was applied to stochastic ODEs by Giles \cite{Giles:2008a,Giles:2008b}. More
recently, it has been used for PDEs with stochastic coefficients \cite{BarthSchwabZollinger:2011,CliffeGilesScheichlTeckentrup:2011}. We
now briefly introduce the MLMC approach in a general context, and derive our MLMC-GMsFEM framework for uncertainty
quantification.

Let $X(\omega)$ be a random variable. We are interested in the efficient computation of
the expected value of $X$, denoted by $\EE[X]$. In our calculations, $X$ is a function of
the permeability field $k$, e.g., the solution to \eqref{eqn:model} evaluated at measurement
points. To compute an approximation to $\EE[X]$, a standard approach is the
Monte Carlo (MC) method. Specifically, one first generates a number $M$ of independent realizations of the
random variable $X$, denoted by $\{X^m\}_{m=1}^M$, and then approximates the expected value
$\EE[X]$ by the arithmetic mean 
$$
E_M(X) := \frac{1}{M} \sum_{m=1}^{M} X^m.
$$
Now we define the Monte Carlo integration error $e_M(X)$ by
\begin{equation*}
  e_M(X)= \EE[X] - E_M(X).
\end{equation*}
Then the central limit theorem asserts that for large $M$, the Monte Carlo integration error
\begin{equation}\label{eqn:MLerror}
  e_M(X)\sim \mathrm{Var}[X]^{1/2}M^{-1/2}\nu,
\end{equation}
where $\nu$ is a standard normal random variable, and $\mathrm{Var}[X]$ is the variance of
$X$. Hence the error $e_M(X)$ in Monte Carlo integration is of order $O(M^{-1/2})$ with a
constant depending only on the variance $\mathrm{Var}[X]$ of the integrand $X$ \cite{RobertCasella:2004}.

In this work, we are interested in MLMC methods. The idea is to compute the quantity of
interest $X=X_L$ using the information on several different levels. Here we couple the MLMC with the
GMsFEM, where the level is identified with the size of the online space. We assume that $L$ is the level
of interest, and computing many realizations at this level is very expensive. Hence we introduce levels
smaller than $L$, namely $L-1,\dots,1$, and assume that the lower the level is, the cheaper the
computation of  $X_l$ is, and the less accurate $X_l$ is with respect to $X_L$. By setting $X_0=0$,
we decompose $X_L$ into
\begin{equation*}
  X_L = \sum_{l=1}^L \left( X_l - X_{l-1}\right).
\end{equation*}
The standard MC approach works with $M$ realizations of the random variable $X_L$ at the level of
interest $L$. In contrast, within the MLMC approach, we work with $M_l$ realizations of $X_l$ at
each level $l$, with $M_1 \geq M_2 \geq \dots \geq M_L$. We write
\begin{equation*}
\EE \left[ X_L \right] = \sum_{l=1}^L \EE \left[ X_l - X_{l-1} \right],
\end{equation*}
and next approximate $\EE \left[ X_l - X_{l-1} \right]$ by an empirical mean:
\begin{equation}\label{eqn:mlmc-l}
  \EE \left[ X_l - X_{l-1} \right]\approx E_{M_l}(X_l - X_{l-1}) = \frac{1}{M_l} \sum_{m=1}^{M_l} \left( X_{l}^m -  X_{l-1}^m \right),
\end{equation}
where $X_l^{m}$ is the $m$th realization of the random variable $X$ computed at the level $l$ (note
that we have $M_l$ copies of $X_l$ and $X_{l-1}$, since $M_l \leq M_{l-1}$). Then the MLMC
approach approximates $\EE[X_L]$ by
\begin{equation}\label{eqn:approx_mlmc}
   E^L(X_L) := \sum_{l=1}^L E_{M_l}\left( X_l - X_{l-1}\right).
\end{equation}
We note that the realizations of $X_l$ used with those of $X_{l-1}$ to evaluate $E_{M_l}
\left( X_l - X_{l-1}\right)$ do not have to be independent of the realizations of $X_l$
used with those of $X_{l+1}$ to evaluate $E_{M_{l+1}}\left(X_{l+1}-X_l\right)$.
In our context, the permeability field samples used for computing $E_{M_l}(X_l-X_{l-1})$ and
$E_{M_{l+1}}(X_{l+1}-X_l)$ do not need to be independent.

We would like to mention that the MLMC can also be interpreted as a multilevel control variate,
following \cite{Speight:2009}. Specifically, suppose that $X=X_L$ on the level $L$ is the quantity
of interest. According to the error estimate \eqref{eqn:MLerror}, the error is proportional to the
product of $M_L^{-1/2}$ and the variance $\mathrm{Var}[X_L]$ of $X_L$. Let $X_{L-1}$
be a cheaper pointwise approximation, e.g., the finite element approximation on a coarser grid,
to $X_L$ with known expected value. Then it is natural
to use $X_{L-1}$ as the control variate to $X_L$ \cite{RobertCasella:2004} and to approximate the expected value $\EE[X_L]$ by
\begin{equation*}
  \begin{aligned}
    \EE[X_L] &= \EE[X_L-X_{L-1}] + \EE[X_{L-1}]\\
      & \approx E_{M_L}(X_L-X_{L-1}) + \EE[X_{L-1}].
    \end{aligned}
\end{equation*}
Here we approximate the expected value $\EE[X_L-X_{L-1}]$ by a Monte Carlo estimate, which,
according to the error estimate \eqref{eqn:MLerror}, will have a small error, if the approximations $X_L$
and $X_{L-1}$ are close to each other. More generally, with a proper choice of weights, the latter condition can be relaxed
to high correlation. In practice, the expected value $\EE[X_{L-1}]$ may be still nontrivial to evaluate.
In the spirit of classical multilevel methods, we can further
approximate the expected value $\EE[X_{L-1}]$ by
\begin{equation*}
  \EE[X_{L-1}] \approx E_{M_{L-1}}(X_{L-1}-X_{L-2}) + \EE[X_{L-2}],
\end{equation*}
where $X_{L-2}$ is a cheap approximation to $X_{L-1}$. By applying this idea recursively,
one arrives at the MLMC estimate as described in \eqref{eqn:approx_mlmc}.

Now we can give the outline of the MLMC-GMsFEM framework, cf. Algorithm \ref{alg:mlmc}. Here,
 the offline space is fixed and preprocessed. The level of the samples
is determined by the size of the online multiscale basis functions. The larger the online
multiscale space $V_\mathrm{on}$ is, the higher the solution resolution is, but the more expensive the
computation is; the smaller the online multiscale space $V_\mathrm{on}$ is, the cheaper the computation
is, but the lower the solution resolution is. The MLMC approach as described above provides
a framework for elegantly combining the hierarchy of approximations from the GMsFEM, and
leveraging the expensive computations on level $L$ to those lower level approximations. In addition,
we note that the samples $\{k^m\}_{m=1}^{M_l}$ used
in the Monte Carlo estimate $ E_{M_{l}}(X_{l}-X_{l-1})$ are identical on every two consecutive
levels, i.e., the permeability samples $\{k^m\}_{m=1}^{M_l}$ used in the Monte Carlo estimates on two consecutive levels are nested.
\begin{algorithm}
\centering
\caption{MLMC-GMsFEM}\label{alg:mlmc}
\begin{itemize}
   \item[1.]  Offline computations
   \begin{itemize}
      \item Construct the snapshot space. 
      \item Construct a low-dimensional offline space by model reduction.
   \end{itemize}
   \item[2.] Multi-level online computations for estimating an expectation at level $l$, $1 \leq l \leq L$.
   \begin{itemize}
      \item Generate $M_l$ realizations of the permeability $\{k_l^{m}\}_{m=1}^{M_l}$ (from
      $\{k_{l-1}^m\}_{m=1}^{M_{l-1}}$).
      \item For each realization $k_l^{m}$, compute online multiscale basis functions.
      \item Solve the coarse-grid problem for $X_{l}^m$.
      \item Calculate the arithmetic mean $E_{M_l}(X_l-X_{l-1})$ by \eqref{eqn:mlmc-l}.
   \end{itemize}
   \item[3.] Output the MLMC approximation $E^L(X)$ by \eqref{eqn:approx_mlmc}.
\end{itemize}
\end{algorithm}

\subsection{Cost analysis}

In the following, we are interested in the root mean square errors
\begin{equation*}
  e_{MLMC}(X_L) = \sqrt{\EE[\|\EE[X_L]-E^L(X_L)\|^2]},
\end{equation*}
\begin{equation*}
  e_{MC}(X_L) = \sqrt{\EE[\|\EE[X_L]-E_{M_L}(X_L)\|^2]},
\end{equation*}
for the MLMC estimate $E^L(X_L)$ and the MC estimate $E_{M_L}(X_L)$, respectively,
with an appropriate norm depending on the quantity of interest (e.g., the absolute value for any
entry of the permeability coefficient, and the $L^2$-norm of the solution). For the error estimation,
we will use the fact that for any random variable $X$ and any norm, $\EE[\|\EE[X]-E_M(X)\|^2]$
defines a norm on the error $\EE[X]-E_M(X)$, and further, there holds the relation
$$
\EE[\|\EE[X]-E_M(X)\|^2]=\frac{1}{M}\EE[\|X-\EE[X]\|^2].
$$

In the analysis, we will be dealing with solutions at different scales. In the MLMC framework,
we denote the scale hierarchy by $H_1\geq H_2\geq \cdots\geq H_L$. The number of realizations
used at the level $l$ for the scale $H_l$ is denoted by $M_l$. We take
\begin{equation*}
  M_1\geq M_2\geq \cdots M_L.
\end{equation*}

For the MLMC approach, the error reads
\begin{equation*}
  \begin{aligned}
    e_{MLMC}(X_L) & = \sqrt{\EE[\|\EE[X_L]-E^L(X_L)\|^2]}\\
     & = \sqrt{\EE[(\EE[\sum_{l=1}^L(X_l-X_{l-1})]-\sum_{l=1}^LE_{M_l}(X_l-X_{l-1}))^2]}\\
     & = \sqrt{\EE[(\sum_{l=1}^L(\EE-E_{M_l})(X_l-X_{l-1}))^2]}\\
     & \leq\sum_{l=1}^L\sqrt{\EE[((\EE-E_{M_l})(X_l-X_{l-1}))^2]}\\
     & \leq\sum_{l=1}^L\frac{1}{\sqrt{M_l}}\sqrt{\EE[(X_l-X_{l-1}-\EE(X_l-X_{l-1}))^2]}
  \end{aligned}
\end{equation*}
where the second last line follows from the triangle inequality for norms, and the last line follows
from \eqref{eqn:MLerror}.
Next we rewrite $X_l-X_{l-1}=(X_l-X)+(X-X_{l-1})$, and since $M_l\leq M_{l-1}$, we deduce
\begin{equation*}
  \begin{aligned}
    e_{MLMC}(X_L) &\leq \sum_{l=1}^L \frac{1}{\sqrt{M_l}}\left(\sqrt{\EE[(X_l-X-\EE(X_l-X))^2]}
        +\sqrt{\EE[(X_{l-1}-X-\EE(X_{l-1}-X))^2]}\right)\\
        & =  \frac{1}{\sqrt{M_L}}\sqrt{\EE[(X_L-X-\EE(X_L-X))^2]} +\frac{1}{\sqrt{M_1}}\sqrt{\EE[X^2]} \\
        &\qquad + \sum_{l=1}^{L-1}\left(\frac{1}{\sqrt{M_{l+1}}}+\frac{1}{\sqrt{M_l}}\right)\sqrt{\EE[((X_l-X)-E(X_l-X))^2]}\\
        & \leq \sum_{l=1}^L\frac{2}{\sqrt{M_{l+1}}}\sqrt{\EE[(X_l-X)^2]}+\frac{1}{\sqrt{M_1}}\sqrt{\EE(X^2)}\\
        &\leq  \sum_{l=1}^L \frac{2}{\sqrt{M_{l+1}}}\delta_l + \frac{1}{\sqrt{M_1}}\sqrt{\EE[X^2]},
  \end{aligned}
\end{equation*}
where the second last line follows from the inequality $\sqrt{\EE[((X_l-X)-E(X_l-X))^2]}\leq
\sqrt{\EE[(X_l-X)^2]}$, a direct consequence of the bias-variance decomposition, and the assumption $M_{L+1}\leq M_L$.
Here we denote $\sqrt{\EE[(X_l-X)^2]}$ as $\delta_l$. As mentioned in section \ref{subsec:mlmcbg}, the lower the
level $l$ is, the less accurate the approximation $X_l$ is with respect to $X_L$, hence we will
have $\delta_1 > \delta_2 > \dots > \delta_L$. To equate the error terms we choose
\begin{equation*}
  M_l = M \left\{
  \begin{array}{l l}
    \left( \frac{1}{\delta_L} \right)^2 \EE[X^2],  & \quad l = 1, \\
   \left(  \frac{\delta_{l-1}}{\delta_L} \right)^2, & \quad 2 \le l \le L+1,
  \end{array} \right.
\end{equation*}
where $M>1$ is a fixed positive integer. Then we end up with
\begin{equation*}
  e_{MLMC}(X_L) \leq \frac{(2L+1)}{M}\delta_L. 
\end{equation*}
In principle, for a prescribed error bound $\epsilon$ such that $e_{MLMC}(X_L)\leq\epsilon$, one can deduce from
the formula the proper choice of the number $N$ of samples, for any given level $L$.

\section{Multilevel Markov chain Monte Carlo}\label{sec:mlmcmc}
One of the most popular and versatile methods for numerically exploring posterior distributions arising from
the Bayesian formulation is the Markov chain Monte Carlo (MCMC) method. The basic idea is to construct
a Markov chain with the target distribution as its stationary distribution. However, the sampling step
remains very challenging in high-dimensional spaces. One powerful idea of improving the sampling
efficiency is preconditioning, first illustrated in \cite{ChristenFox:2005,EfendievHouLuo:2006},
and more recently extended in \cite{BalLangmore:2013}. In the latter work, some theoretical
properties, e.g., asymptotic confidence interval, of a multistage version of the two-level
algorithm \cite{EfendievHouLuo:2006} are also established. 

\subsection{MLMCMC with GMsFEM}

The standard Metropolis-Hastings algorithm generates samples from the posterior distribution
$\pi(k)=p(k|F_{obs})$, cf. Algorithm \ref{alg:mcmc}. Here $\mathcal{U}(0,1)$ is the uniform distribution
over the interval $(0,1)$. As was described in Section \ref{sec:prelim}, the permeability field $k$ is determined by
the parameters $\theta$ and $\tau$. Hence, given the current sample $k^m$, parameterized by its parameters
$\theta^m$ and $\tau^m$, one can generate the proposal $k$ by generating the proposal for
$\theta$ and $\tau$ first, i.e., draw $\theta$
from distribution $q_{\theta}(\theta | \theta^m) $ and $\tau$ from distribution
$q_{\tau}(\tau | \tau^m) $ (in view of the independence between $\theta$ and $\tau$), for some
proposal distributions $q_\theta(\theta|\theta^m)$ and $q_\tau(\tau|\tau^m)$, and then
form the proposal for the entire permeability field $k$.

\begin{algorithm}
    \caption{Metropolis-Hastings MCMC} \label{alg:mcmc}
    \begin{algorithmic}[1]
      \STATE Specify $k_0$ and $M$.
      \FOR {$m=0:M$}
      \STATE Generate the entire permeability field proposal $k$ from $q(k|k^m)$.
      \STATE Compute the acceptance probability $\gamma(k^m)$ by \eqref{eqn:mcmc-accept}.
      \STATE Draw $u\sim \mathcal{U}(0,1)$.
      \IF {$\gamma(k^m,k)\leq u$}
         \STATE $k^{m+1}=k$.
      \ELSE
         \STATE $k^{m+1}=k^m$.
      \ENDIF
      \ENDFOR
    \end{algorithmic}
\end{algorithm}

The transition kernel $K_r(k^m,k)$ of the Markov chain generated by Algorithm \ref{alg:mcmc} is given by
\begin{equation*}
  K_r(k^m,k) = \gamma(k^m,k) q(k|k^m) + \delta_{k^m}(k) \left(1- \int \gamma(k^m,k) q(k|k^m) dk \right),
\end{equation*}
where $q(k|k^m)$ denotes the proposal distribution and $\gamma(k^m,k)$ denotes the
acceptance probability for the proposal $k$ defined by
\begin{equation}\label{eqn:mcmc-accept}
     \gamma(k^m,k) = \min \left\{ 1, \frac{q(k^m|k)\pi(k)}{q(k|k^m)\pi(k^m)} \right\}.
\end{equation}

Now we integrate the multilevel idea with the Metropolis-Hastings algorithm and the GMsFEM.
Like before, we start with the telescopic sum
\begin{equation*}
  \begin{aligned}
    \EE_{\pi_L}[F_L] & = \int F_L (x)\pi_L(x) dx \\
      & =\int F_0(x)\pi_0(x) dx +\sum_{l=1}^L \int(F_l(x)\pi_l(x)-F_{l-1}(x)\pi_{l-1}(x))dx,
  \end{aligned}
\end{equation*}
where $\pi_l$ denotes the approximated target distribution at level $l$, and $\pi_0$ is our initial level.
We note that after the initial level each expectation involves two measures,
$\pi_l$ and $\pi_{l-1}$, which is different from
the case of the MLMC (see \cite{ketelson2013}). Therefore, we rewrite the integration using a product measure as
\begin{equation*}
  \begin{aligned}
   \int(F_l(x)\pi_l(x)-F_{l-1}(x)\pi_{l-1}(x))dx &= \int F_l(x)\pi_l(x) dx - \int F_{l-1}(y)\pi_{l-1}(y)dy\\
     & = \int\int(F_l(x)-F_{l-1}(y))\pi_l(x)\pi_{l-1}(y)dxdy\\
     & = \EE_{\pi_l,\pi_{l-1}}[F_l(x)-F_{l-1}(y)].
  \end{aligned}
\end{equation*}
Therefore, we have
\begin{equation} \label{eqn:mlmcmcsum}
    \EE_{\pi^L}[F_L] =\EE_{\pi_0}[F_0]+ \sum_{l=1}^L \EE_{\pi_l,\pi_{l-1}}[F_l-F_{l-1}].
\end{equation}
The idea of our multilevel method is to estimate each term of the right hand side of equation \eqref{eqn:mlmcmcsum} independently.
In particular we can estimate each term in \eqref{eqn:mlmcmcsum} by an MCMC estimator. The first term
$\EE_{\pi_0}[F_0]$ can be estimated using the standard MCMC estimator in Algorithm \ref{alg:mcmc}.
We estimate the expectation $\EE_{\pi_l,\pi_{l-1}}[F_l(x)-F_{l-1}(y)]$ by the sample mean
\begin{equation} \label{eqn:estimatorlevel}
   \EE_{\pi_l,\pi_{l-1}}[F_l(x)-F_{l-1}(y)] \approx \frac{1}{M_l}\sum_{m=1}^{M_l}(F_l(x_l^m)-F_{l-1}(y_l^m)),
\end{equation}
where the samples $\{(y_l^m,x_l^m)\}_{m=1}^{M_l}$ are drawn from the product measure $\pi_{l-1}(y)\otimes\pi_l(x)$.
Next we describe an efficient preconditioned MCMC method for generating samples from the product measure $\pi_{l-1}(y)\otimes\pi_l(x)$, extending
our earlier work \cite{EfendievHouLuo:2006}.

Here we introduce a multilevel MCMC algorithm by adapting the proposal distribution $q(k|k_m)$ to
the target distribution $\pi(k)$ using the GMsFEM with different sizes of the online space which we call different
levels, cf. Algorithm \ref{alg:mlmcmc}. The process modifies the proposal distribution $q(k|k_m)$
by incorporating the online coarse-scale information. Let $F_l(k)$ be the pressure/production computed by solving
online coarse problem at level $l$ for a given $k$. The target distribution $\pi(k)$ is approximated on
level $l$ by $\pi_l(k)$, with $\pi(k)\equiv\pi_L(k)$. Here we have
\begin{equation} \label{eqn:pi}
   \pi_l(k) \propto \exp \left( - \frac{ || F_{obs} - F_l(k)|| ^2}{2\sigma_l^2} \right) \times p(k).
\end{equation}

In the algorithm we still keep the same offline space for each level. From level $0$ to level $L$,
we increase the size of the online space as we go to a higher level, which means for any levels
$l$, $l+1\le L$, samples of level $l$ are cheaper to generate than that of level $l+1$. This idea underlies
the cost reduction using the multilevel estimator. Hence the
posterior distribution for coarser levels $\pi_l, l = 0,\dots, L-1$ do not have to model the
measured data as faithfully as $\pi_L$, which in particular implies that by choosing suitable value of $\sigma_l^2$ it
is easier to match the result $F_l(k)$ with the observed data. We denote the number of samples at
level $l$ by $M_l$, where we will have $M_0 \le \dots \le M_L$. As was discussed above, our quantity
of interest can be approximated by the telescopic sum \eqref{eqn:mlmcmcsum}. We denote the estimator of
$\EE_{\pi_0}[F_0]$ at the initial level by $\widehat{F}_0$. Then by the MCMC estimator we have
\begin{equation*}
\widehat{F}_0 = \frac{1}{M_0}\sum_{m =1}^{M_0} F_0(x_0^m).
\end{equation*}
Here $x^m_0$ denotes the samples we accepted on the initial level (after discarding the samples
at the burn-in period). Similarly we denote the estimator of the differences $\EE_{\pi_l,
\pi_{l-1}}[F_l(x)-F_{l-1}(y)]$ by $\widehat{Q}_l$. Then with \eqref{eqn:estimatorlevel} we have
$$
\widehat{Q}_l = \frac{1}{M_l}\sum_{m=1}^{M_l}(F_l(x_l^m)-F_{l-1}(y_l^m)),
$$
where the samples
$\{(y_l^m,x_l^m)\}_{m=1}^{M_l}$ are drawn from the product measure $\pi_{l-1}(y)\otimes\pi_l(x)$.
Finally denote the estimator of $\EE_{\pi_L}[F_L]$ or our full MLMCMC estimator by $\widehat{F}_L$,
then the quantity of interest $E_{\pi_L}[F_L]$ is approximated by
\begin{equation} \label{eqn:mlmcmc_est}
   \widehat{F}_L =\widehat{F}_0 + \sum_{l=1}^L  \widehat{Q}_l.
\end{equation}

\begin{algorithm}
  \centering
  \caption{Multilevel Metropolis-Hastings MCMC}\label{alg:mlmcmc}
  \begin{algorithmic}[1]
	     \STATE Given $k_m$, draw a trial proposal $k$ from distribution $q(k|k_1^m) =
           q_0(k|k_1^m)$
	     \STATE Compute the acceptance probability
         \[
           \rho_1(k_1^m,k) = \min \left\{ 1, \frac{q_0(k_1^m|k)\pi_1(k)}{q_0(k|k_1^m)\pi_1(k_1^m)} \right\}
         \]
         \STATE $u\sim \mathcal{U}(0,1)$
         \IF {$u<\rho_1(k_1^m,k)$}
           \STATE $k_1^{m+1}=k$ (at the initial level)
         \ELSE
            \STATE $k_1^{m+1}=k_1^m$ (at the initial level)
         \ENDIF
         \FOR {$l=1:L-1$}
           \IF {$k$ is accepted at level $l$}
              \STATE Form the proposal distribution $q_{l}$ (on the $l+1$th level) by
                 \begin{equation*}
                    q_{l}(k|k_{l+1}^m) = \rho_l (k_{l+1}^m,k)q_{l-1}(k|k_{l+1}^m) +\delta_{k_{l+1}^m} (1 - \int \rho_l (k_{l+1}^m,k)q_{l-1}(k|k_{l+1}^m) dk_{l+1}^m )
                 \end{equation*}
              \STATE Compute the acceptance probability
                 \begin{equation*}
                    \rho_{l+1}(k_{l+1}^m,k) = \min \left \{ 1, \frac {q_l(k_{l+1}^m|k) \pi_{l+1}(k)}{q_l(k|k_{l+1}^m) \pi_{l+1}(k_{l+1}^m)} \right \} = \min \left \{ 1, \frac {\pi_l(k_{l+1}^m) \pi_{l+1}(k)}{\pi_l(k) \pi_{l+1}(k_{l+1}^m)} \right \}
                 \end{equation*}
              \STATE $u\sim \mathcal{U}(0,1)$
              \IF {$u<\rho_{l+1}(k_{l+1}^m,k)$}
                  \STATE $k_{l+1}^{m+1}=k$ and go to next level (if $l=L-1$, accept $k$ and set $k_L^{m+1}=k$).
              \ELSE
                  \STATE $k_{l+1}^{m+1}=k_{l+1}^m$, and break
              \ENDIF
           \ENDIF
         \ENDFOR
   \end{algorithmic}
\end{algorithm}

\subsection{Convergence analysis}
In this part, we briefly analyze the convergence property of the multilevel MCMC algorithm, cf.
Algorithm \ref{alg:mlmcmc}. Specifically, we discuss the detailed balance relation and
the ergodicity of the Markov chain, following the general convergence theory in \cite{RobertCasella:2004}.

To this end, we denote
\begin{equation}
  \begin{aligned}
  \bigE^l &= \{ k: \pi_l (k) > 0 \},\quad  l = 1,2,\ldots, L, \\
  \bigD & =  \{ k: q_{l-1}(k|k_l^m) >0 \quad\text{for some} \,\, k_l^m \in \bigE^l \}. 
 \end{aligned}
\end{equation}
The set $\bigE^l$ is the support of the distributions $\pi_l(k)$ at level $l$. The set $\bigE^L$
is the support of the target distribution $\pi(k) = \pi_L(k)$ at the finest level. The set
$\bigD$ is the set of all the proposals which can be generated by the proposal distribution
$q_{l-1}(k| k_l^m)$. To sample from $\pi(k)$ correctly, it is necessary that $\bigE^L \subseteq
\bigE^{L-1} \subseteq \ldots \subseteq \bigE^1 \subseteq \bigD$ (up to a set of zero measure).
Otherwise, if one of these conditions is not true, say, $\bigE^{l+1} \not\subseteq \bigE^l$,
then there will exist a subset $A \subset (\bigE^{l+1} \setminus \bigE^l)$ such that
\begin{equation*}
  \pi_{l+1}(A) = \int_A \pi_{l+1}(k) dk >0 \quad \text{and} \quad  \pi_l(A) = \int_A \pi_l(k) dk = 0,
\end{equation*}
which means no element of $A$ can pass the level $l$ and $A$ will never be visited by the Markov chain
$\{ k_{l+1}^m\}$. Thus, the distribution at level $l+1$, i.e., $\pi_{l+1}(k)$ is not sampled properly.

For most practical proposal distributions $q_{l-1}(k|k_l^m)$, such as random walk samplers, the condition $\bigE^1,\ldots,
\bigE^L \subseteq \bigD$ is naturally satisfied. To show the inclusion $\bigE^{l+1} \subseteq \bigE^l$ for any level
$l$, notice that if the precision parameters $\sigma_{l+1}$ and $\sigma_l$ are chosen to be relatively
small, then $\pi_{l+1}(k)$ and $ \pi_l(k)$ are very close to zero for most proposals. From the numerical
point of view, the proposal $k$ is very unlikely to be accepted if $\pi_{l+1}(k)$ and $ \pi_l(k)$ are
close to zero. Consequently the support of the distributions should be interpreted as
\begin{equation*}
   \bigE^{l+1} = \{ k : \pi_{l+1} (k) > \delta \}\quad\mbox{ and }\quad  \bigE^{l} = \{ k: \pi_{l} (k) > \delta \}
\end{equation*}
where $\delta$ is a small positive number. If $k \in \bigE^{l+1}$, then $\pi_{l+1}(k) > \delta$ and
$\| F_{obs} - F_{l+1}(k)\|^2/ 2\sigma^2_{l+1} $ is not very large. To make $k\in \bigE^l$,  $\| F_{obs}
- F_l(k)\|^2/ 2\sigma^2_{l} $ should not be very large either. If $\| F_{obs} - F_l(k)\| $ is bounded
by $\| F_{obs} - F_{l+1}(k) \| $ up to a multiplicative constant, then the condition $\bigE^{l+1}
\subseteq \bigE^l$ can be satisfied by choosing the parameter $\sigma_l$ properly. For our model,
the coarser level quantity is indeed bounded by the fine level quantity.
Thus, the condition $\bigE^L \subseteq \bigE^{L-1} \subseteq \ldots \subseteq \bigE^1 \subseteq \bigD$ is satisfied.

Let
\begin{equation}
  Q_l(k_l^m,k) = \rho_l (k_l^m,k)q_{l-1}(k |k_l^m) +\delta_{k_l^m} (1 - \int \rho_l (k_l^m,k)q_{l-1}(k | k_l^m) dk_l^m) \label{trans_ker}
\end{equation}
denote the transition kernel of the Markov chain at level $l$.
As in a regular MCMC method, we can show that $Q_l(k_l^m,k)$ satisfies the detailed balance condition
at level $l$, i.e.,
\begin{equation}
  \pi_{l}(k_l^m)Q_l(k_l^m,k) = \pi_{l}(k) Q_l(k,k_l^m) \label{balance_cond}
\end{equation}
for any $k,k_l^m \in \bigE^l$. In fact, the equality \eqref{balance_cond} is obviously true when $k = k_l^m$.
If $k \neq k_l^m$, then $Q_l(k_l^m,k) = \rho_{l} (k_l^m,k)q_{l-1}(k | k_l^m)$, we have
\begin{equation*}
  \begin{aligned}
   \pi_{l}(k_l^m)Q_l(k_l^m,k) &= \pi_{l}(k_l^m)\rho_{l} (k_l^m,k)q_{l-1}(k | k_l^m)\\
    & = \min \left( \pi_{l}(k_l^m)q_{l-1}(k | k_l^m), \pi_{l}(k)q_{l-1}(k_l^m |k) \right)\\
    &= \min \left( \frac{\pi_{l}(k_l^m)q_{l-1}(k | k_l^m)}{\pi_{l}(k)q_{l-1}(k_l^m| k)} ,1 \right)   \pi_{l}(k)q_{l-1}(k_l^m| k)\\
    &= \rho_{l} (k,k_l^m) \pi_{l}(k) q_{l-1}(k_l^m| k) = \pi_{l}(k)Q_{l}(k,k_l^m).
  \end{aligned}
\end{equation*}
So the detailed balance condition \eqref{balance_cond} is always satisfied. Using
\eqref{balance_cond} we can easily show that $\pi(A) = \int Q_l(k,A) dk$ for any $A \in \mathcal{B}(\bigE^l)$,
where $\mathcal{B}(\bigE^l)$ denotes all measurable subsets of $\bigE^l$. Thus, $\pi_l(k)$ is indeed
the stationary distribution of the transition kernel $Q_l(k_l^m,k)$.

In a regular MCMC method, cf. Algorithm \ref{alg:mcmc}, the proposal $q(k|k^m) $ is usually chosen to satisfy $q(k|k^m) >0$
for any $(k^m,k)\in \bigE \times \bigE$, which guarantees that the resulting MCMC chain is
irreducible. Similarly the irreducibility holds for multilevel MCMC at each level $l$ if
$q_{l-1}(k|k_l^m) >0$ for any $(k_l^m,k)\in \bigE^l \times \bigE^l$. We already have  $\bigE^L
\subseteq \bigE^{L-1} \subseteq \ldots \subseteq \bigE^1 $  holds, which means $\rho_l(k_l^m,k)
>0$, and also for common choices of the proposal distribution, we have $q_{l-1}(k|k_l^m)$
positive, which guarantees the irreducibility of the chain at each level.

To prove the convergence  of the distribution, we need to show that the chain is aperiodic. Recall
that a simple sufficient condition for aperiodicity is that the transition kernel $Q(k^m,\{k^m\}) > 0 $
for some $k^m \in \bigE$. In other words, the event $\{ k^{m+1} = k^m\}$ happens with a positive
probability. For our multilevel MCMC at finest level $l$, consider the transition kernel \eqref{trans_ker}, we have
\begin{eqnarray*}
 Q_l(k_l^m,\{k_l^m\}) &=& 1 - \int_{k \neq k_l^m}  \rho_l (k_l^m,k)q_{l-1}(k | k_l^m) dk_l^m  \\ &=&
         	1 - \int_{k \neq k_l^m}  \rho_l (k_l^m,k) \rho_{l-1}(k_l^m ,k) \dots \rho_1(k_l^m, k) q_{0}(k|k_l^m) dk_l^m
\end{eqnarray*}
Hence $Q_l(k_l^m,\{k_l^m\})  \equiv 0$ requires $\rho_s (k_s^m,k) = 1$ for $s = 1, \ldots, l$, for almost all
$k \in \bigD$. which means that all the proposals generated by $q_{0}(k_1^m,k)$ are correct samples
for distributions at all levels. In this case it does not make sense to use the MCMC method since we can sample
directly from $q(k|k^m)$. Thus in practice we can always safely assume that the chain generated by the
multilevel MCMC is aperiodic. As a result the Markov chain generated by MLMCMC converges.

In Algorithm \ref{alg:mlmcmc}, the specific proposal distribution $q_l$ can be computed easily and at no additional cost,
as we can simplify the acceptance probability for level $l+1$ to
\begin{equation}\label{eqn:accept_mlmcmc}
  \rho_{l+1}(k_{l+1}^m,k) = \min \left \{ 1, \frac {\pi_l(k_{l+1}^m) \pi_{l+1}(k)}{\pi_l(k) \pi_{l+1}(k_{l+1}^m)} \right \}.
\end{equation}
This is true when $k_{l+1}^m = k$, so we will demonstrate this for the case $k_{l+1}^m \neq k$. In this case,
$q_{l}(k|k_{l+1}^m) = \rho_{l}(k_{l+1}^m,k)q_{l-1}(k|k_{l+1}^m) $,  then,
\begin{equation*}
    \rho_{l+1}(k_{l+1}^m,k) = \min \left \{ 1, \frac {q_l(k_{l+1}^m|k) \pi_{l+1}(k)}{q_l(k|k_{l+1}^m) \pi_{l+1}(k_{l+1}^m)} \right \}
           =  \min \left \{ 1, \frac {\rho_l(k,k_{l+1}^m) q_{l-1}(k_{l+1}^m|k) \pi_{l+1}(k)}
           {\rho_l(k_{l+1}^m,k) q_{l-1}(k|k_{l+1}^m) \pi_{l+1}(k_{l+1}^m)} \right \}.
\end{equation*}
Assume for simplicity $ q_{l-1}(k_{l+1}^m|k) \pi_{l}(k) > q_{l-1}(k|k_{l+1}^m) \pi_{l}(k_{l+1}^m)$, then
$\rho_l(k_{l+1}^m,k) = 1$ and $\rho_l(k,k_{l+1}^m) = \frac{ q_{l-1}(k|k_{l+1}^m) \pi_{l}(k_{l+1}^m)}{ q_{l-1}
(k_{l+1}^m|k) \pi_{l}(k)}$. Using these relations we obtain the desired formula \eqref{eqn:accept_mlmcmc}.
Similarly, in the case of $ q_{l-1}(k_{l+1}^m|k) \pi_{l}(k) < q_{l-1}(k|k_{l+1}^m) \pi_{l}(k_{l+1}^m)$, then $\rho_l(k,k_{l+1}^m) = 1$ and
$\rho_l(k_{l+1}^m,k) = \frac{ q_{l-1}(k_{l+1}^m|k) \pi_{l}(k)}{ q_{l-1}(k|k_{l+1}^m) \pi_{l}(k_{l+1}^m)}$. With these
relations we also deduce that \eqref{eqn:accept_mlmcmc} holds.

\section{Numerical results}
\label{sec:numer}

In our numerical examples, we consider permeability fields described by two-point correlation functions, and use
Karhunen-Lo\`eve expansion (KLE) to parameterize the permeability
fields. Then we apply the MLMC and MLMCMC with the GMsFEM algorithms described earlier.
First, we briefly recall the permeability parametrization, and then we present numerical results.

\subsection{Permeability parameterization}

To obtain a permeability field in terms of an optimal $L^2$ basis, we use the KLE \cite{Loeve:1977}.
For our numerical tests, we truncate the expansion and represent the permeability matrix by a
finite number of random parameters. We consider the random field $Y( {x},\omega)=\log[k( {x},\omega)]$, where
$\omega$ represents randomness. We assume a zero mean $\EE[Y( {x},\omega)]=0$, with a known covariance
operator $R({x},{y})=\EE\left[Y( {x})Y( {y})\right]$. Then we expand the random field $Y(x,\omega)$ as
\[
  Y( {x},\omega)=\sum_{k=1}^{\infty} Y_k(\omega) \Phi_k( {x}),\qquad
\]
with
\begin{equation*}
Y_k(\omega)=\int_\Omega Y(x,\omega)\Phi_k(x)dx.
\end{equation*}
The functions $\{\Phi_k(x)\}$ are eigenvectors of the covariance operator $R(x,y)$, and form
a complete orthonormal basis in $L^2(\Omega)$, i.e.,
\begin{equation}\label{eig}
\int_{\Omega} R( {x}, {y})\Phi_k( {y})d {y}=\lambda_k\Phi_k( {x}),
\qquad k=1,2,\ldots,
\end{equation}
where $\lambda_k=\EE[Y_k^2]>0$. We note that $\EE[Y_iY_j]=0$ for
all $i\neq j$. By denoting $\eta_k=Y_k/\sqrt{\lambda_k}$
(whence $\EE[\eta_k]=0$ and $\EE[\eta_i\eta_j]=\delta_{ij}$), we have
\begin{equation}\label{KLE}
   Y( {x},\omega)=\sum_{k=1}^{\infty} \sqrt{\lambda_k}\eta_k(\omega)\Phi_k( {x}),
\end{equation}
where $\Phi_k$ and $\lambda_k$ satisfy \eqref{eig}. The randomness is represented
by the scalar random variables $\eta_k$. After discretizing the domain $\Omega$ into
a rectangular mesh, we truncate the KLE \eqref{KLE} to a finite number of terms.
In other words, we keep only the leading-order terms (quantified by the magnitude
of $\lambda_k$), and  capture most of the energy of the stochastic process $Y(x,
\omega)$. For an $N$-term KLE approximation
\begin{equation*}
  Y_N=\sum_{k=1}^{N}\sqrt{\lambda_k}\eta_k\Phi_k,
\end{equation*}
the energy ratio of the approximation is defined by
\begin{equation*}
  e(N):=\frac{E\|Y_N\|^2}{E\|Y\|^2}=\frac{\sum_{k=1}^N\lambda_k}{\sum_{k=1}^{\infty}\lambda_k}.
\end{equation*}
If the eigenvalues $\{\lambda_k\}$ decay very fast, then the truncated KLE with the
first few terms would be a good approximation of the stochastic process $Y(x,\omega)$ in the $L^2$ sense.

In our examples, the permeability field $k$ is assumed to follow a log-normal distribution
with a known spatial covariance, with the correlation function $R(x,y)$ given by
\begin{equation}\label{normal}
   R(x,y)=\sigma^2\exp{\Bigl(-\frac{|x_1-y_1|^2}{2l_1^2}-\frac{|x_2-y_2|^2}{2l_2^2}\Bigr)}.
\end{equation}
where $l_1$ and $l_2$ are the correlation lengths in $x_1$- and $x_2$-direction,
respectively, and $\sigma^2=\EE[Y^2]$ is a constant that determines the variation
of the permeability field.

\subsection{MLMC}

In our simulations, we evaluate the performance of the MLMC method on computing the
expected values of our quantity of interest $F$. In particular, we consider the
stationary, single-phase flow model \eqref{eqn:model} on the unit square domain
$\Omega = (0,1)^2$ with $f \equiv 1$ and linear boundary conditions. The forward
problem is solved with the GMsFEM, and the fine grid and coarse grid are chosen to be $
50 \times 50$ and $5 \times 5$, respectively. The quantity of interest $F$ for this
set of simulations is the fine scale pressure field. We consider the following
two Gaussian covariance functions:
\begin{itemize}
 \item Isotropic Gaussian field with correlation lengths $l_1= l_2=0.1$ and a stochastic dimension 5;
 \item Anisotropic Gaussian field with correlation lengths $l_1=0.1$ and $l_2=0.05$, and a stochastic dimension 5.
\end{itemize}
In both cases, we use the variance $\sigma^2 = 2$, and keep $N=5$ terms in the final KL expansion
where the $\eta_k$ coefficients are drawn from a normal distribution with zero mean and unit variance.

\begin{figure}
        \centering
        \begin{subfigure}[b]{1\textwidth}
                \centering
                \includegraphics[trim = .5cm 1cm .5cm .5cm, clip=true,width=\textwidth]{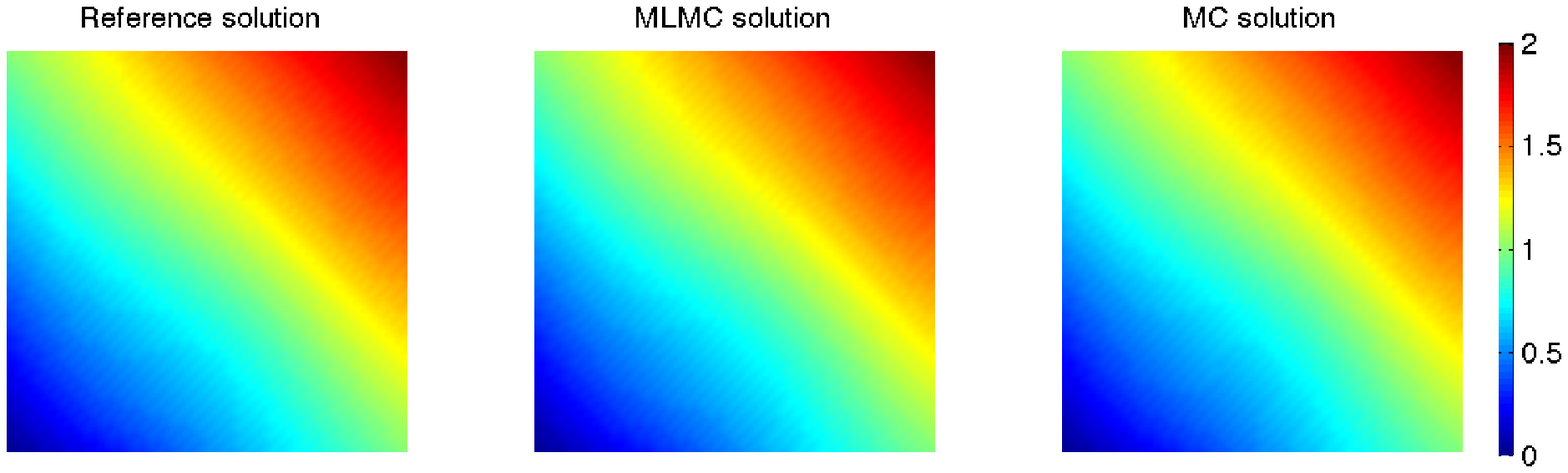}
                \caption{Isotropic Gaussian}
                \label{fig:isoMLMC}
        \end{subfigure}

        \begin{subfigure}[b]{1\textwidth}
                \centering
                \includegraphics[trim = .5cm 1cm .5cm .5cm, clip=true,width=\textwidth]{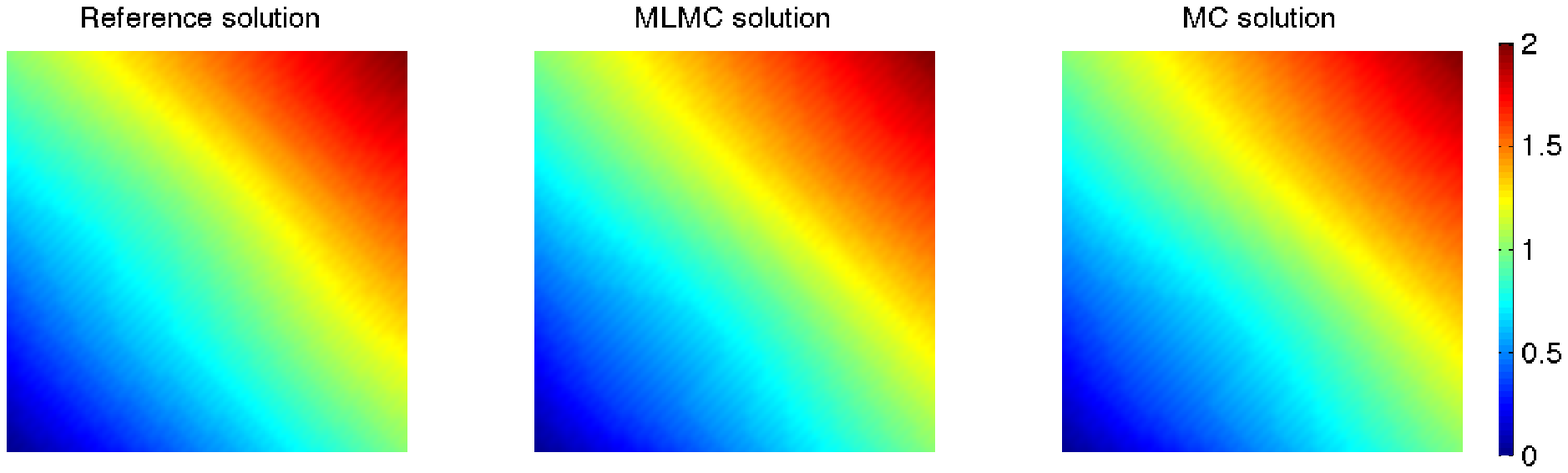}
                \caption{Anisotropic Gaussian}
                \label{fig:anisoMLMC}
        \end{subfigure}

\caption{Pressure field solutions for different methods and correlation lengths}\label{fig:MLMC}
\end{figure}

We denote by $F_l$ the fine scale pressure field at level $l$ in MLMC. The level of our interest is $L=3$.
As stated in Algorithm \ref{alg:mlmc}, we generate $M_l$ realizations at level $l$ of the permeability field,
solve the model problems by choosing $N_l$ eigenvalues to generate the online space in the GMsFEM, and
compute the MLMC approximation of $\EE[F_L]$ by \eqref{eqn:approx_mlmc}. We compare MLMC with the standard MC at the
level $L$ of interest with the same amount of cost. Hence we choose
\begin{equation*}
  \widehat{M} = \frac{\sum^L_{l=1} N^2_l M_l}{N^2_L}
\end{equation*}
as the number of samples in the standard MC algorithm. We use the arithmetic mean of
$M_{ref}$ samples of the pressure field as reference and compute the relative $L^2$-errors
\begin{equation*}
  (e^{rel}_{MLMC})[F_L] = \frac{\|\EE^{ref}_{M_{ref}}[F_L] - \EE^L[F_L] \| _{L^2(D)}}{\|\EE^{ref}_{M_{ref}}[F_L] \| _{L^2(D)}}
\end{equation*}
\begin{equation*}
   (e^{rel}_{MC})[F_L]
= \frac{\|\EE^{ref}_{M_{ref}}[F_L] - \EE^{MC}_{\widehat{M}}[F_L] \| _{L^2(D)}}{\|\EE^{ref}_{M_{ref}}[F_L] \| _{L^2(D)}}.
\end{equation*}
For the simulations use $N_1 = 4$, $N_2 = 8$, and $N_3 = 16$ eigenfunctions for the online
space construction. We respectively set the number of samples at each level to be $M_1 = 128$,
$M_2 = 32$, and $M_3 = 8$, and equate the computational costs for the MLMC and MC relative
error comparisons. With this choice of realizations for MLMC, we use $\widehat{M} = 20$
permeability realizations for the standard MC forward simulations. The parameters we have
used and the respective relative errors are summarized in Table \ref{MLMCtable}. Figure
\ref{fig:MLMC} illustrates expected  pressure fields for different correlation lengths and
different methods (MLMC and MC). For both covariances, we observe that the MLMC approach yields
errors which are about 1.5 times smaller than those resulting from the MC approach. We note
that the gain is larger for the isotropic case than for the anisotropic case.

\begin{table}
\begin{center}
 \caption{Parameters and errors for the estimates by MLMC vs. MC} \label{MLMCtable}
\begin{tabular}{lcc}
 \hline
                                     & Isotropic Gaussian & Anisotropic Gaussian\\
 \hline
 $(N_1, N_2, N_3)$  &   $(4,8,16)$              &    $(4,8,16)$               \\
 $(M_1, M_2, M_3)$ &   $(128, 32, 8) $      &      $(128, 32, 8) $             \\
 $\widehat{N} $                    &    16                           &           16        \\
 $\widehat{M}$                    &    24                           &         24           \\
 $ M_{MCref}$            &   5000                       &         5000            \\
 $e^{rel}_{MLMC}$             &    0.0431                  &       0.0653            \\
 $e^{rel}_{MC}$                 &    0.0802                  &       0.0952            \\
 $e^{rel}_{MC} / e^{rel}_{MLMC}$      &   1.86    &   1.45    \\ \hline
 \end{tabular}
  \end{center}
 \end{table}

\subsection{MLMCMC}

\begin{figure}
        \centering
        \includegraphics[trim = 0cm .5cm .5cm .5cm, clip=true,width=0.33\textwidth]{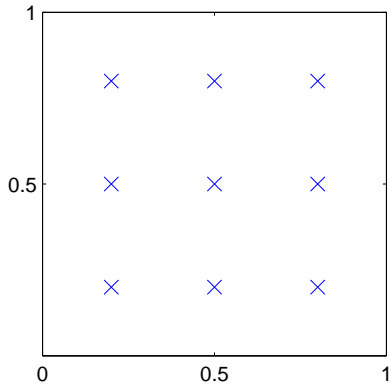}
 \caption{The points where pressure is evaluated}\label{fig:pressurept}
 \end{figure}

In our MLMCMC experiment we also consider the model problem ~\eqref{eqn:model} on $\Omega = (0,1)^2$
with $f \equiv 1$ and linear boundary conditions. The prior permeability distribution $p(k)$ is also
parameterized by KLE as above. The ``observed" data $F_{obs}$ is obtained by generating a reference
permeability field, solving the forward problem with the GMsFEM, and evaluating the pressure at
nine points away from the boundary. The locations of the reference pressures are shown in Figure \ref{fig:pressurept}.

Our proposal distribution is a random walker sampler in which the proposal distribution depends on
the previous value of the permeability field and is given by $q(k|k_n) = k_n + \delta \epsilon_n $
where $\epsilon_n$ is a random perturbation with mean zero and unit variance, and $\delta$ is a
step size. The random perturbations are imposed on the $\eta_k$ coefficients in the KL expansion.

We consider two examples, one with isotropic Gaussian field of correlation length $l_1= l_2=0.1$,
the other with anisotropic Gaussian field of correlation lengths $l_1= 0.05$,  $l_2=0.1$.
For both examples we use $\delta = 0.2$ in the random walk sampler.
We again use the level $L = 3$, and for each level $l$
we take the same number of KLE terms, $N = 5$ for the tests.
For the GMsFEM, we take the number of eigenvalues
to generate the online space at each level as $N_1 = 4, N_2 = 8, N_3 = 16$.
We take our quantities of interest $F$ as the pressure values at the same nine
points and use them in order to compute the acceptance probabilities as shown in Algorithm \ref{alg:mlmcmc}.

\begin{figure}
   \centering
   \includegraphics[trim = .5cm .5cm .5cm .5cm, clip=true,width=0.67\textwidth]{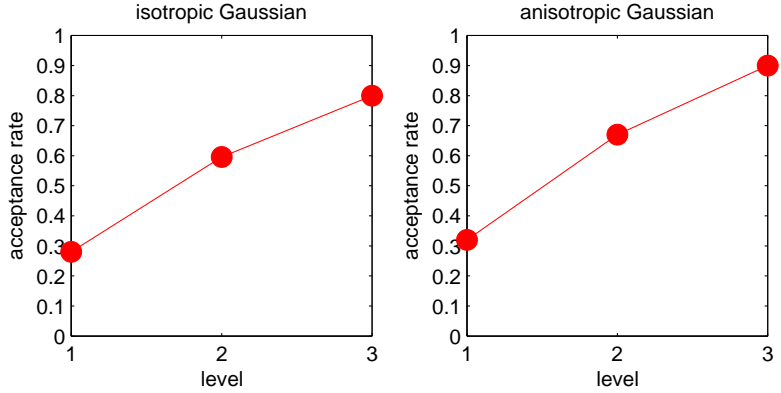}
   \caption{Acceptance rate of multilevel sampler with both isotropic and anisotropic trials}\label{fig:accpt_rate}
 \end{figure}
\begin{figure}
  \centering
  \includegraphics[trim = .5cm .5cm .5cm .5cm, clip=true,width=0.68\textwidth]{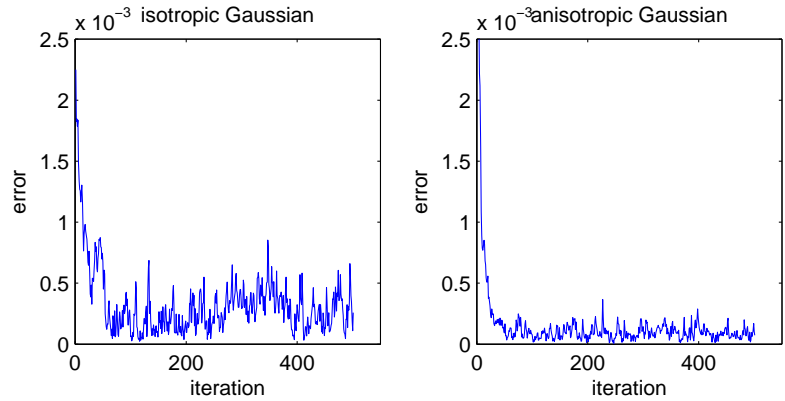}
  \caption{Plots of iteration vs. error with both isotropic and anisotropic tirals}\label{fig:err_iter}
\end{figure}
\begin{figure}
        \centering
                \includegraphics[trim = .5cm 1cm .5cm .5cm, clip=true,width=\textwidth]{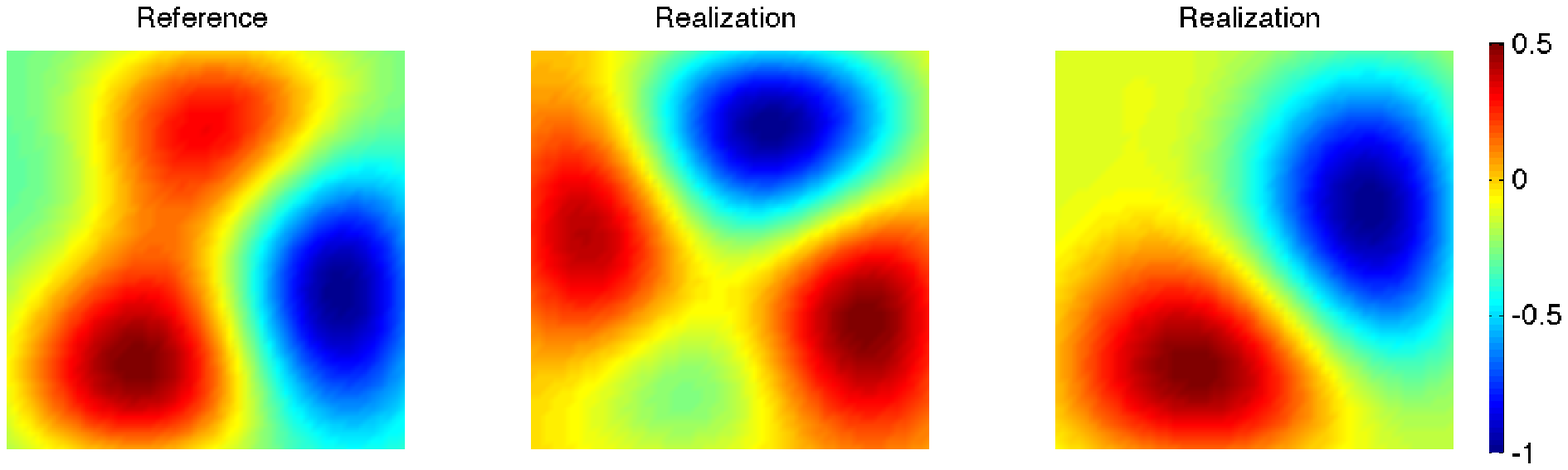}

                \includegraphics[trim = .5cm 1cm .5cm .5cm, clip=true,width=\textwidth]{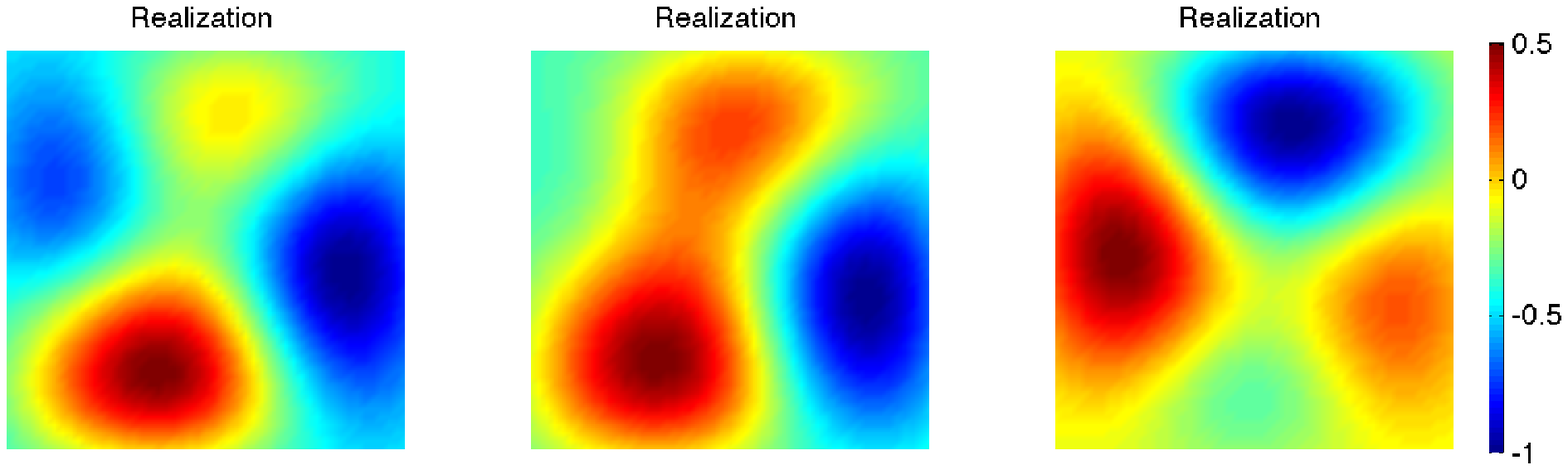}

\caption{Isotropic MLMCMC accepted realization}\label{fig:isoMLMCMC}
\end{figure}

For the MLMCMC examples, we run Algorithm \ref{alg:mlmcmc} until $P_4 = 1000$ total samples
pass the final level of acceptance. We note that 300 initial accepted samples are discarded as burn-in.
The acceptance rates of the multilevel sampler are shown in Fig. \ref{fig:accpt_rate}. To
compute the acceptance rates, we assume that $P_1$, $P_2$, and $P_3$ samples are proposed for
respective levels $L_1$, $L_2$, and $L_3$. Then, the rate at the $l$-th level is the ratio
$P_{l+1} / P_l$. Most notably, the results in Fig. \ref{fig:accpt_rate} show that
the acceptance rate increases as $l$ increases. In particular, for more
expensive (larger) levels, we observe that it is much more probable that a
proposed sample will be accepted. This is an advantage of the multilevel method,
due to the fact that less proposals are wasted on more expensive computations.  We also show
a set of plots in Fig. \ref{fig:err_iter} that illustrate the errors $ E_k = \|F_{obs} -
F_k \|$, cf. \eqref{eqn:pi}, of the accepted samples on the finest level. In Fig.
\ref{fig:isoMLMCMC} we plot some of the accepted permeability realizations that have
passed all levels of computation. We note that the general shapes of the accepted the
fields do not necessarily match that of the reference field, reinforcing the notion
that the problem is ill-posed due to the fact that a variety of proposals may explain
the reference data equally well.

\section{Concluding remarks}
\label{sec:concl}
In this paper we propose a robust framework for the uncertainty quantification of quantities of
interest for high-contrast single-phase flow problems. The procedure combines the generalized
multiscale finite element method (GMsFEM) and multilevel Monte Carlo (MLMC) methods. Within this
context, GMsFEM provides a hierarchy of approximations at varying levels of accuracy and
computational cost, and MLMC offers an efficient way to estimate quantities of interest using
samples on respective levels. The number of basis functions in the online GMsFEM stage may be
readily and adaptively modified in order to adjust the computational costs and requisite accuracy, and efficiently
generate samples at different levels. In particular, it is cheap to generate
samples through smaller dimensional online spaces with less accuracy, and it is expensive to
generate samples through larger dimensional spaces with a higher level of accuracy. As such,
a suitable choice of the number of samples at different levels allows us to
leverage the expensive computations at finer levels toward the coarse grid, while retaining the
accuracy of the final estimates on the output quantities of interest. We additionally describe a
multilevel Markov chain Monte Carlo (MLMCMC) inverse modeling technique, which sequentially screens
the proposal with different levels of GMsFEM approximations. In particular, the method reduces the
number of evaluations that are required at finer levels, while combining the samples at varying
levels to arrive at an accurate estimate. A number of numerical examples are presented in order
to illustrate the efficiency and accuracy of multilevel methods as compared to standard Monte Carlo estimates.
The analysis and examples of the proposed methodology offer a seamless integration between the
flexibility of the GMsFEM online space construction along with the multilevel features of MLMC methods.

\bibliographystyle{abbrv}
\bibliography{mlmcmc}

\end{document}